\documentclass[11pt,reqno]{amsart}
\usepackage{amsfonts,amssymb,amsmath}

\DeclareMathOperator{\bmax}{\mathbf{max}}
\DeclareMathOperator{\trace}{\mathrm{trace}}
\DeclareMathOperator{\conv}{\mathrm{conv}}
\DeclareMathOperator{\ex}{\mathrm{ex}}
\DeclareMathOperator{\KG}{\mathsf{KG}}

\newtheorem{theorem}{Theorem}

\newtheorem{proposition}[theorem]{Proposition}

\newtheorem{example}[theorem]{Example}

\numberwithin{equation}{section} \numberwithin{theorem}{section}

\begin{document}

\title{Pattern Recognition on Oriented Matroids: Three--Tope Committees}

\author{Andrey O. Matveev}
\address{Data-Center Co., RU-620034, P.O.~Box~5, Ekaterinburg,
Russian~Federation} \email{aomatveev@hotmail.com aomatveev@dc.ru}


\begin{abstract}
A three-tope committee $\mathcal{K}^{\ast}$ for a simple oriented matroid~$\mathcal{M}$ is a~$3$-subset of its maximal covectors such that every positive halfspace of~$\mathcal{M}$ contains at least two topes from $\mathcal{K}^{\ast}$. We consider three-tope committees as the~vertex sets of triangles in graphs associated with the~topes and enumerate them making use of the properties of the poset of convex subsets of the ground set of $\mathcal{M}$.
\end{abstract}

\maketitle

\pagestyle{myheadings}

\markboth{PATTERN RECOGNITION ON ORIENTED MATROIDS}{A.O.~MATVEEV}

\thispagestyle{empty}

\tableofcontents

\section{Introduction}

Throughout this note, which is a sequel to~\cite{AM1,AM2}, $\mathcal{M}:=(E,\mathcal{T})$ denotes a~simple oriented matroid that is not acyclic, on ground set $E$, with set of topes~$\mathcal{T}$; we say that~$\mathcal{M}$ is {\em simple\/} if it has no loops, parallel or {\sl antiparallel\/} elements. For an element~\mbox{$e\in E$}, we denote by $\mathcal{T}_e^+:=\{T\in\mathcal{T}:\ T(e)=+\}$ the {\em positive halfspace\/} of~$\mathcal{M}$, corresponding to $e$.

A {\em tope committee\/} $\mathcal{K}^{\ast}$ {\em for $\mathcal{M}$} is a subset of its maximal covectors such that
$|\{K\in\mathcal{K}^{\ast}:\ K(e)=+\}|>|\{K\in\mathcal{K}^{\ast}:\ K(e)=-\}|$ or, in other words, $|\mathcal{T}_e^+\cap\mathcal{K}^{\ast}|>\tfrac{1}{2}|\mathcal{K}^{\ast}|$, for each element $e\in E$.
This construction can serve as an abstract analogue of building blocks of decision rules in contradictory problems of pattern recognition~\cite{AM1}.

A comprehensive survey of (in)feasibility studies is given in~\cite{Ch}.

When the oriented matroid $\mathcal{M}$ interprets in the pattern recognition problem as (a reorientation of) an abstract {\em training set}, the three-tope committee becomes the~most preferred approximation to the notion of {\em solution\/} of an~infeasible system of constraints, and we restrict our attention here to three-tope committees.

General three-tope committees are considered and enumerated in Section~\ref{sec:1}, while Section~\ref{sec:2} is devoted to the committees (that are the~most important constructions for recognition goals) whose topes have inclusion-maximal positive parts.

We use the following notation: $|\cdot|$ and $\#$ denote the cardinality of a set and the number of sets in a family, respectively. If $\mathcal{F}$ is a set of sign vectors, then $-\mathcal{F}:=\{-F:\ F\in\mathcal{F}\}$. We let $F^+$ and $F^-$ denote the positive and negative parts of a sign vector $F$, respectively.

For a bounded poset $P$ we denote by $P^{\mathrm{c}}$ the set of its coatoms. The~singleton has no atoms and coatoms. For a chain $x\lessdot\ y$ of unit length we set $\{x,y\}^{\mathrm c}:=\{x\}$. See, e.g.,~\cite[\S{}IV.4.B]{A} on a link between the number of coatoms in $P$ and incidence functions on $P$. $[x,z]:=\{y\in P:\ x\leq y\leq z\}$ is a closed interval in $P$, and $\mu_P(\cdot,\cdot)$ denotes the M\"{o}bius function on $P$.

See, e.g.,~\cite[Chapter~9]{BLSWZ},\cite{E1,E2,ERW,LV,S} on {\em convexity\/} in oriented matroids. We adopt from~\cite{ERW} the following terminology and statements: A~subset \mbox{$A\subset E$} is {\em acyclic\/} if the restriction of $\mathcal{M}$ to $A$ is an acyclic oriented matroid. The~{\em convex hull\/}~$\conv(A)$ of an acyclic set $A$ is the set $B\supseteq A$ such that for every covector~$F$ and for every element $b\in B$, the implication $F(a)=+$, $\forall a\in A$ $\Longrightarrow$ $F(b)=+$ holds. An acyclic set $A$ is {\em convex\/} if $\conv(A)=A$; its subset of {\em extreme points\/} $\ex(A)$ is defined by $\ex(A):=\bigl\{a\in A:\ a\not\in\conv(A-\{a\})\bigr\}$. For any acyclic set $A$ it holds $\conv(\ex(A))=\conv(A)$. A convex set $A\subset E$ is {\em free\/} if $\ex(A)=\conv(A)=A$. The~meet-semilattice $L_{\conv}(\mathcal{M})$ is defined as the~family of convex subsets of~$E$ ordered by inclusion.

We regard the least element~$\hat{0}$ of $L_{\conv}(\mathcal{M})$ as the empty subset of~$E$. If a subset $H\subseteq E$ is not acyclic then we set $\conv(H):=E$ and~$\ex(E):=E$.

For a simple graph $\pmb{G}$ we let $\mathfrak{V}(\pmb{G})$ and $\mathfrak{E}(\pmb{G})$ denote its vertex set and the~family of edges, respectively.

If $\mathcal{F}$ is a set family then the {\em Kneser graph\/} $\KG(\mathcal{F})$ of $\mathcal{F}$
is the graph with $\mathfrak{V}\bigl(\KG(\mathcal{F})\bigr):=\mathcal{F}$; if $F',F''\in\mathcal{F}$ then $\{F',F''\}$ \mbox{$\in\mathfrak{E}\bigl(\KG(\mathcal{F})\bigr)$} iff $|F'\cap F''|=0$, see~\cite[\S{}3.3]{MBZ}.

\section{Three-Tope (Anti-)Committees}
\label{sec:1}

A {\em tope anti-committee\/} $\mathcal{A}^{\ast}$ for $\mathcal{M}$ is a subset of its maximal covectors such that the set $-\mathcal{A}^{\ast}$ is a tope committee for $\mathcal{M}$.
We denote by $\mathbf{K}_k^{\ast}(\mathcal{M})$ and~$\mathbf{A}_k^{\ast}(\mathcal{M})$ the families of all tope committees and anti-committees of cardinality $k$ for $\mathcal{M}$, respectively. Recall that due to axiomatic symmetry~(L1)~\cite[\S4.1.1]{BLSWZ} we have
\begin{equation}
\label{eq:3}
\mathcal{T}=-\mathcal{T}
\end{equation}
and thus it follows from the definition that
\begin{equation*}
\#\mathbf{K}_k^{\ast}(\mathcal{M})=\#\mathbf{A}_k^{\ast}(\mathcal{M})=
\#\mathbf{K}_{|\mathcal{T}|-k}^{\ast}(\mathcal{M})=\#\mathbf{A}_{|\mathcal{T}|-k}^{\ast}(\mathcal{M})
\end{equation*}
because for any tope committee $\mathcal{K}^{\ast}$ for $\mathcal{M}$ its complement $\mathcal{T}-\mathcal{K}^{\ast}$ is an~anti-committee.

Since a three-tope anti-committee meets each positive halfspace in at most one tope, we have
\begin{equation}
\label{eq:1}
\#\mathbf{A}_3^{\ast}(\mathcal{M})=
\sum_{\substack{S_1\subset E,\ S_2\subset E,\ S_3\subset E:\\ |S_1|,|S_2|,|S_3|>0,\\ |S_1\cap S_2|=|S_1\cap S_3|=|S_2\cap S_3|=0}}\ \prod_{k\in\{1,2,3\}}\ \ \Bigl|\bigcap_{s\in S_k}\mathcal{T}_s^+ - \bigcup_{s\in E-S_k}\mathcal{T}_s^+\Bigr|\ ;
\end{equation}
note that $\bigcap_{s\in S_k}\mathcal{T}_s^+ - \bigcup_{s\in E - S_k}\mathcal{T}_s^+=
\bigcap_{s\in S_k}\mathcal{T}_s^+ \cap \Bigl(-\bigcap_{s\in E - S_k}\mathcal{T}_s^+\Bigr)$.

On the right-hand side of~(\ref{eq:1}) the (extreme points of) convex sets $S_1$, $S_2$ and $S_3$ only contribute to the total sum, therefore we turn to the poset $L_{\conv}(\mathcal{M})$ of convex subsets of the ground set of $\mathcal{M}$.

Let $\widehat{L}_{\conv}(\mathcal{M})$ be the poset $L_{\conv}(\mathcal{M})$ augmented by a greatest element~$\hat{1}$; it is convenient to set $\hat{1}:=E$. For brevity, we will write $\widehat{L}$ instead of~$\widehat{L}_{\conv}(\mathcal{M})$.

If we let $\mathcal{T}_B^+:=\bigcap_{b\in B}\mathcal{T}_b^+$ denote the intersection of the positive halfspaces corresponding to the elements of a subset $B\subseteq E$, then expression~(\ref{eq:1}) reformulates in the following way:
\begin{multline}
\label{eq:8}
\#\mathbf{K}_3^{\ast}(\mathcal{M})=\#\mathbf{A}_3^{\ast}(\mathcal{M})\\=
\sum_{\substack{\{A_1,A_2,A_3\}\subset \widehat{L}-\{\hat{0},\hat{1}\}:\\
A_1\wedge A_2=A_1\wedge A_3=A_2\wedge A_3=\hat{0}}}\ \prod_{k\in\{1,2,3\}}\ \ \Bigl|\mathcal{T}_{\ex(A_k)}^+ \cap \left(-\mathcal{T}_{\ex(\conv(E-A_k))}^+\right)\Bigr|\ .
\end{multline}
The related quantity
\begin{equation*}
8\binom{|\mathcal{T}|/2}{3}+\sum_{A\in\widehat{L}-\{\hat{0},\hat{1}\}}\mu_{\widehat{L}}(\hat{0},A)
\binom{|\mathcal{T}_{\ex(A)}^+|}{3}=
8\binom{|\mathcal{T}|/2}{3}+\sum_{\substack{A\in\widehat{L}-\{\hat{0},\hat{1}\}:\\ \text{$A$ free}}}(-1)^{|A|}
\binom{|\mathcal{T}_A^+|}{3}
\end{equation*}
is the number of three-tope sets containing no opposites and meeting every positive halfspace of $\mathcal{M}$. Recall that if $A\in\widehat{L}-\{\hat{0},\hat{1}\}$ then \mbox{$\mu_{\widehat{L}}(\hat{0},A)=(-1)^{|A|}$} whenever the convex set $A$ is free and, as a consequence, $[\hat{0},A]$ is order-isomorphic to~the~Boolean lattice of rank $|A|$; otherwise, \mbox{$\mu_{\widehat{L}}(\hat{0},A)=0$}, see~\cite{ERW}.

Let $\mathbf{G}:=\mathbf{G}(\mathcal{M})$ be a graph with the vertex set $\mathfrak{V}(\mathbf{G}):=\mathcal{T}$; a pair $\{T',T''\}\subset\mathcal{T}$ by definition belongs to the edge family $\mathfrak{E}(\mathbf{G})$ of $\mathbf{G}$ iff no~positive halfspace of $\mathcal{M}$ contains this pair, that is, $|(T')^+\cap (T'')^+|=0$ or, in other words, $\{T',T''\}$ is a $1$-dimensional {\em missing face\/} of the abstract simplicial {\em complex $\Delta_{\text{\rm acyclic}}(\mathcal{M})$ of acyclic subsets\/} of $E$. Thus, $\mathbf{G}$ is isomorphic to the Kneser graph $\KG\bigl(\{T^+:\ T\in\mathcal{T}\}\bigr)$ of the family of the positive parts of topes of $\mathcal{M}$.

Let $\boldsymbol{\Gamma}:=\boldsymbol{\Gamma}(\mathcal{M})$ be the graph defined by
\begin{gather*}
\mathfrak{V}(\boldsymbol{\Gamma}):=\mathcal{T}\ ,\\
\{T',T''\}\in\mathfrak{E}(\boldsymbol{\Gamma})\ \ \Longleftrightarrow\ \ (T')^+\cup (T'')^+=E\ .
\end{gather*}
The vertex set of any odd cycle in $\boldsymbol{\Gamma}$ is a tope committee for $\mathcal{M}$ \cite[\S5]{AM1}. Since
$\{T',T''\}\in\mathfrak{E}(\mathbf{G})$ iff $(T')^-\cup (T'')^-=E$, for any pair of topes $\{T',T''\}$,
the~four graphs $\boldsymbol{\Gamma}$, $\mathbf{G}$, $\KG\bigl(\{T^+:\ T\in\mathcal{T}\}\bigr)$ and $\KG\bigl(\{T^-:\ T\in\mathcal{T}\}\bigr)$ are all isomorphic, due to symmetry~(\ref{eq:3}).
For example, the mapping $\mathfrak{V}(\mathbf{G})\to\mathfrak{V}(\boldsymbol{\Gamma})$, $T\mapsto-T$, is an isomorphism between the graphs $\mathbf{G}$ and $\boldsymbol{\Gamma}$. Since a~subset $\{T',T'',T'''\}\subset\mathcal{T}$ is a three-tope anti-committee for $\mathcal{M}$ iff it is the~vertex set of a triangle in $\mathbf{G}$ (or, in other terms, this subset is a \mbox{$2$-dimensional} face of the {\em independence complex \/} of the graph whose edges are the $1$-dimensional faces of~$\Delta_{\text{\rm acyclic}}(\mathcal{M})$), the~three-tope committees for $\mathcal{M}$ are precisely the~vertex sets of the triangles in $\boldsymbol{\Gamma}$. From the poset-theoretic point of view, the~family~$\mathbf{K}_3^{\ast}(\mathcal{M})$ is regarded in~\cite[Proposition~4.1, Theorem~5.1]{AM2} as antichains in posets associated with the topes.

We have
\begin{align}
\label{eq:5}
\#\mathfrak{E}(\boldsymbol{\Gamma})&=
\binom{|\mathcal{T}|}{2}+\sum_{\substack{A\in\widehat{L}-\{\hat{0},\hat{1}\}:\\ \text{$A$ free}}}(-1)^{|A|}
\binom{|\mathcal{T}_A^+|}{2}
\intertext{and}
\label{eq:6}
|\mathfrak{V}(\boldsymbol{\Gamma})|:=|\mathcal{T}|&=\sum_{\substack{A\in\widehat{L}-\{\hat{0},\hat{1}\}:\\ \text{$A$ free}}}(-1)^{|A|-1}
\bigl|\mathcal{T}_A^+\bigr|\ ;
\end{align}
recall that in fact the most efficient tool of tope enumeration is the Las~Vergnas--Zaslavsky formula, see~\cite[\S4.6]{BLSWZ}.

Now let $\mathbf{R}$ be a graph, with the vertex set $\mathfrak{V}(\mathbf{R}):=\{1,2,\ldots,|\mathcal{T}|\}$, isomorphic to either of the graphs $\boldsymbol{\Gamma}$, $\mathbf{G}$, $\KG\bigl(\{T^+:\ T\in\mathcal{T}\}\bigr)$ and $\KG\bigl(\{T^-:\ T\in\mathcal{T}\}\bigr)$. If we let~$\mathrm{A}$ and $\mathcal{N}(i)$ denote the adjacency matrix of $\mathbf{R}$ and the neighborhood of the vertex $i$ in $\mathbf{R}$, respectively, then well-known observations of graph theory imply, for example, that
$\#\mathbf{K}_3^{\ast}(\mathcal{M})=
\tfrac{1}{6}\trace(\mathrm{A}^3)$ and
\begin{equation}
\label{eq:4}
\#\mathbf{K}_3^{\ast}(\mathcal{M})=\frac{1}{3}
\sum_{\{i,j\}\in\mathfrak{E}(\mathbf{R})}|\mathcal{N}(i)\cap\mathcal{N}(j)|
\end{equation}
because the quantity $\#\mathbf{K}_3^{\ast}(\mathcal{M})=\#\mathbf{A}_3^{\ast}(\mathcal{M})$ is the number of triangles in $\mathbf{R}$.

\section{Committees of Cardinality Three whose Topes Have Maximal Positive Parts}
\label{sec:2}

Let $\bmax^+(\mathcal{T})$ denote the set of topes of the oriented matroid $\mathcal{M}$ whose positive parts are maximal with respect to inclusion; these parts are the~inclu\-sion-maximal convex subsets of the ground set~$E$ and thus they are the~co\-atoms of the lattice~$\widehat{L}$, that is, $\{T^+:\ T\in\bmax^+(\mathcal{T})\}=\widehat{L}^{\mathrm c}$. If $\mathcal{M}$ is realizable by a~central arrangement of oriented hyperplanes, then the coatoms of~$\widehat{L}$ are the multi-indices of maximal feasible subsystems of a certain related homogeneous system of strict linear inequalities. The {\em graph of topes with maximal positive parts\/} $\boldsymbol{\Gamma}_{\bmax}^+:=\boldsymbol{\Gamma}_{\bmax}^+(\mathcal{M})$ is the subgraph of~$\boldsymbol{\Gamma}(\mathcal{M})$ with the~vertex set $\bmax^+(\mathcal{T})$, that is,
\begin{gather*}
\mathfrak{V}(\boldsymbol{\Gamma}_{\bmax}^+):=\bmax^+(\mathcal{T})\ ,\\
\{T',T''\}\in\mathfrak{E}(\boldsymbol{\Gamma}_{\bmax}^+)\ \ \Longleftrightarrow\ \ (T')^+\cup(T'')^+=E\ .
\end{gather*}
The graph~$\boldsymbol{\Gamma}_{\bmax}^+$ is connected~\cite[\S5.2]{AM1}.

Recall that if  $T\in\bmax^+(\mathcal{T})$ then symmetry~(\ref{eq:3}) and maximality of the~positive part $T^+\in\widehat{L}^{\mathrm c}$ imply that the set $E-T^+=T^-$ is acyclic and convex; as a consequence, we have
\begin{multline*}
\#\{\mathcal{K}^{\ast}\in\mathbf{K}_3^{\ast}(\mathcal{M}):\ \mathcal{K}^{\ast}\subseteq\bmax^+(\mathcal{T})\} =\#\bigl\{\{D_1,D_2,D_3\}\subset\widehat{L}:\\
\{E-D_1,E-D_2,E-D_3\}\subseteq\widehat{L}^{\mathrm{c}},\
D_1\wedge D_2=D_1\wedge D_3=D_2\wedge D_3=\hat{0}\bigr\}\ ,
\end{multline*}
cf.~(\ref{eq:8}).

The degree of a vertex~$T$ in the graph $\boldsymbol{\Gamma}_{\bmax}^+$ equals the number of coatoms in the interval $[T^-,\hat{1}]$ of the lattice $\widehat{L}$. As a consequence,
the number of edges in $\boldsymbol{\Gamma}_{\bmax}^+$ is
\begin{equation*}
\#\mathfrak{E}(\boldsymbol{\Gamma}_{\bmax}^+)=\frac{1}{2}\sum_{D\in\widehat{L}:\ E-D\in\widehat{L}^{\mathrm{c}}}\left|[D,\hat{1}]^{\mathrm{c}}\right|\ ,
\end{equation*}
and the {\em cyclomatic number\/} of $\boldsymbol{\Gamma}_{\bmax}^+$ equals
\begin{equation*}
1+\frac{1}{2}\sum_{D\in\widehat{L}:\ E-D\in\widehat{L}^{\mathrm{c}}}\left|[D,\hat{1}]^{\mathrm{c}}\right|-
|\widehat{L}^{\mathrm c}|\ .
\end{equation*}

If a tope $T'\in\mathcal{T}$ does not belong to the set $\bmax^+(\mathcal{T})$ then there exists, again thanks to symmetry~(\ref{eq:3}), a tope $T''\in\bmax^+(\mathcal{T})$ such that the~pair~$\{T',T''\}$ is an edge in $\boldsymbol{\Gamma}$; since the graph $\boldsymbol{\Gamma}_{\bmax}^+$ is connected, this implies that the graph $\boldsymbol{\Gamma}$ is connected as well and, according to~(\ref{eq:5}) and~(\ref{eq:6}),
the {\em cyclomatic number\/} of $\boldsymbol{\Gamma}$ equals
\begin{equation*}
1+\binom{|\mathcal{T}|}{2}+\sum_{\substack{A\in\widehat{L}-\{\hat{0},\hat{1}\}:\\ \text{$A$ free}}}(-1)^{|A|}\binom{1+|\mathcal{T}_A^+|}{2}\ .
\end{equation*}

The concluding statement is a direct consequence of expression~(\ref{eq:4}):
\begin{proposition} The number of committees of cardinality three, for the~oriented matroid $\mathcal{M}$, whose topes have inclusion-maximal positive parts is
\begin{equation}
\label{eq:7}
\#\{\mathcal{K}^{\ast}\in\mathbf{K}_3^{\ast}(\mathcal{M}):\ \mathcal{K}^{\ast}\subseteq\bmax^+(\mathcal{T})\}=\frac{1}{3} \sum_{\substack{\{D_1,D_2\}\subset\widehat{L}:\\ \{E-D_1,E-D_2\}\subset\widehat{L}^{\mathrm{c}},\\
D_1\wedge D_2=\hat{0}}}\left|[D_1\vee D_2,\hat{1}]^{\mathrm{c}}\right|\ .
\end{equation}
\end{proposition}

\begin{example} Let $\mathcal{M}:=(E_6,\mathcal{T})$ be the rank\/ $3$ simple oriented matroid on the~ground set $E_6:=\{1,2,\ldots,6\}$, with the set of topes
{\tiny
\begin{equation*}
\mathcal{T}:=
\begin{matrix}
\{&-&-&+&+&+&+\\ &-&-&+&-&+&+\\ &+&-&+&-&+&+\\ &+&-&+&-&+&-&\\
&-&-&+&-&+&-\\ &-&-&+&+&+&-\\ &-&-&+&+&-&+\\ &-&+&+&+&-&+&\\
&-&+&+&+&+&-\\ &-&+&+&-&+&-\\ &+&+&+&-&+&-\\ &+&+&+&-&-&-&\\
&-&+&+&-&-&-\\ &-&+&+&+&-&-\\ &+&-&-&-&+&+\\ &+&-&-&+&+&+&\\
&-&-&-&+&+&+\\ &-&-&-&+&-&+\\ &+&-&-&+&-&+\\ &+&-&-&-&-&+&\\
&+&-&-&-&+&-\\ &+&+&-&-&+&-\\ &+&+&-&-&-&+\\ &+&+&-&+&-&+&\\
&-&+&-&+&-&+\\ &-&+&-&+&-&-\\ &+&+&-&+&-&-\\ &+&+&-&-&-&-&\}\ .
\end{matrix}
\end{equation*}
}
\noindent A realization, by a central arrangement of oriented hyperplanes in $\mathbb{R}^3$, of the~acyclic reorientation ${}_{-\{1,2\}}\mathcal{M}$ of $\mathcal{M}$ is shown in~\cite[Figure~3.1]{AM1}. We have
{\small
\begin{equation*}
\{D\in\widehat{L}:\ E_6-D\in\widehat{L}^{\mathrm c}\}\\=
\bigl\{\{12\},\{15\},\{16\},\{23\},\{24\},\{35\},\{46\}\bigr\}
\end{equation*}
}
and
{\small
\begin{align*}
\left|[\{12\}\vee\{35\},\hat{1}]^{\mathrm c}\right|=\left|[\{1235\},\hat{1}]^{\mathrm c}\right|
=\#\bigl\{\{1235\}\bigr\}&=1\ ,\\
\left|[\{12\}\vee\{46\},\hat{1}]^{\mathrm c}\right|=\left|[\{1246\},\hat{1}]^{\mathrm c}\right|
=\#\bigl\{\{1246\}\bigr\}&=1\ ,\\
\left|[\{15\}\vee\{23\},\hat{1}]^{\mathrm c}\right|=\left|[\{1235\},\hat{1}]^{\mathrm c}\right|
=\#\bigl\{\{1235\}\bigr\}&=1\ ,\\
\left|[\{15\}\vee\{24\},\hat{1}]^{\mathrm c}\right|=\left|\{\hat{1}\}^{\mathrm c}\right|
&=0\ ,\\
\left|[\{15\}\vee\{46\},\hat{1}]^{\mathrm c}\right|=\left|[\{1456\},\hat{1}]^{\mathrm c}\right|
=\#\bigl\{\{1456\}\bigr\}&=1\ ,\\
\left|[\{16\}\vee\{23\},\hat{1}]^{\mathrm c}\right|=\left|\{\hat{1}\}^{\mathrm c}\right|
&=0\ ,
\end{align*}
\begin{align*}
\left|[\{16\}\vee\{24\},\hat{1}]^{\mathrm c}\right|=\left|[\{1246\},\hat{1}]^{\mathrm c}\right|
=\#\bigl\{\{1246\}\bigr\}&=1\ ,\\
\left|[\{16\}\vee\{35\},\hat{1}]^{\mathrm c}\right|=\left|[\{1356\},\hat{1}]^{\mathrm c}\right|
=\#\bigl\{\{1356\}\bigr\}&=1\ ,\\
\left|[\{23\}\vee\{46\},\hat{1}]^{\mathrm c}\right|=\left|[\{2346\},\hat{1}]^{\mathrm c}\right|
=\#\bigl\{\{2346\}\bigr\}&=1\ ,\\
\left|[\{24\}\vee\{35\},\hat{1}]^{\mathrm c}\right|=\left|[\{2345\},\hat{1}]^{\mathrm c}\right|
=\#\bigl\{\{2345\}\bigr\}&=1\ ,\\
\left|[\{35\}\vee\{46\},\hat{1}]^{\mathrm c}\right|=\left|[\{3456\},\hat{1}]^{\mathrm c}\right|
=\#\bigl\{\{3456\}\bigr\}&=1\ .
\end{align*}
}
\noindent Thus, according to~{\rm(\ref{eq:7})}, the family $\{\mathcal{K}^{\ast}\in\mathbf{K}_3^{\ast}(\mathcal{M}):\ \mathcal{K}^{\ast}\subseteq\bmax^+(\mathcal{T})\}$ consists of\/ $\tfrac{9}{3}=3$ committees which can be found as the vertex sets of the~triangles in the graph depicted in~\cite[Figure~5.4]{AM1}.
\end{example}

\end{document}